\documentclass[11pt]{amsart}

\usepackage{latexsym}
\usepackage{amsmath}
\usepackage{amssymb}
\usepackage{epsfig}
\usepackage{amsfonts}

\newcommand\risS[6]{\raisebox{#1pt}[#5pt][#6pt]{\begin{picture}(#4,15)(0,0)
  \put(0,0){\includegraphics[width=#4pt]{#2.eps}} #3
     \end{picture}}}
\newcommand\smf[1]{\risS{-4}{#1}{}{15}{0}{0}}
\newcommand\laf[1]{\begin{picture}(65,50)(0,0)
          \put(0,0){\includegraphics[width=65pt]{#1.eps}}
                   \end{picture}}
\newcommand\baf[1]{\begin{picture}(50,55)(0,0)
          \put(0,0){\includegraphics[width=50pt]{#1.eps}}
                   \end{picture}}
\newcommand\saf[2]{\begin{picture}(#2,55)(0,0)
          \put(0,3){\includegraphics[width=#2pt]{#1.eps}}
                   \end{picture}}
\def\rb#1#2{\raisebox{#1pt}{#2}}
\def\kb#1{[ #1 ]}
\def\wt#1{\widetilde{#1}}
\def\a{\alpha}
\def\b{\beta}
\def\bc{\mathrm{bc}}

\def\wh{\widehat}
\def\wo{\overline}
\def\bu{\bullet}

\def\Ga{\Gamma}

\def\de{\delta}

\def\al{\alpha}
\def\be{\beta}

\def\ve{\varepsilon}

\def\vp{\varphi}

\def\cS{\mathcal S}

\def\cF{\mathcal F}

\def\wt{\widetilde}
\def\<{\langle}
\def\>{\rangle}

\newtheorem{thm}{Theorem}[section]

\newtheorem{cor}[thm]{Corollary}

\newtheorem{defn}[thm]{Definition}

\newtheorem{Example}[thm]{Example}

\begin{document}

\hfill{\it Dedicated to Askold Khovanskii on the occasion of his 60th birthday
       \vspace{10pt}}

\title[The Kauffman bracket and the Bollob\'as-Riordan polynomial]%
    {The Kauffman bracket of virtual links and the Bollob\'as-Riordan polynomial}
\author[Sergei~Chmutov]{Sergei~Chmutov}
\author[Igor~Pak]{Igor~Pak}
\date{}


\begin{abstract}
\noindent
We show that the Kauffman bracket~$\kb{L}$ of a checkerboard
colorable virtual link~$L$ is an evaluation of the
Bollob\'as-Riordan polynomial $R_{G_L}$ of a ribbon graph associated with~$L$.
This result generalizes the celebrated relation between the classical
Kauffman bracket and the Tutte polynomial of planar graphs.
\end{abstract}

\maketitle

\section*{Introduction} \label{s:intro}

The theory of \emph{virtual links} was discovered independently by
L.~Kauffman~\cite{Ka3} and M.~ Goussarov, M.~Polyak, and O.~Viro~\cite{GPV}.
Virtual links are represented by their diagrams which differ from
ordinary knot diagrams by presence of {\it virtual crossings},
which should be understood not as
crossings but rather as defects of our two-dimensional picture. They should
be treated in the same way as the extra crossings appearing in planar pictures of
non-planar graphs. Virtual link diagrams are considered modulo the {\it classical}
Reidemeister moves
$$\risS{-18}{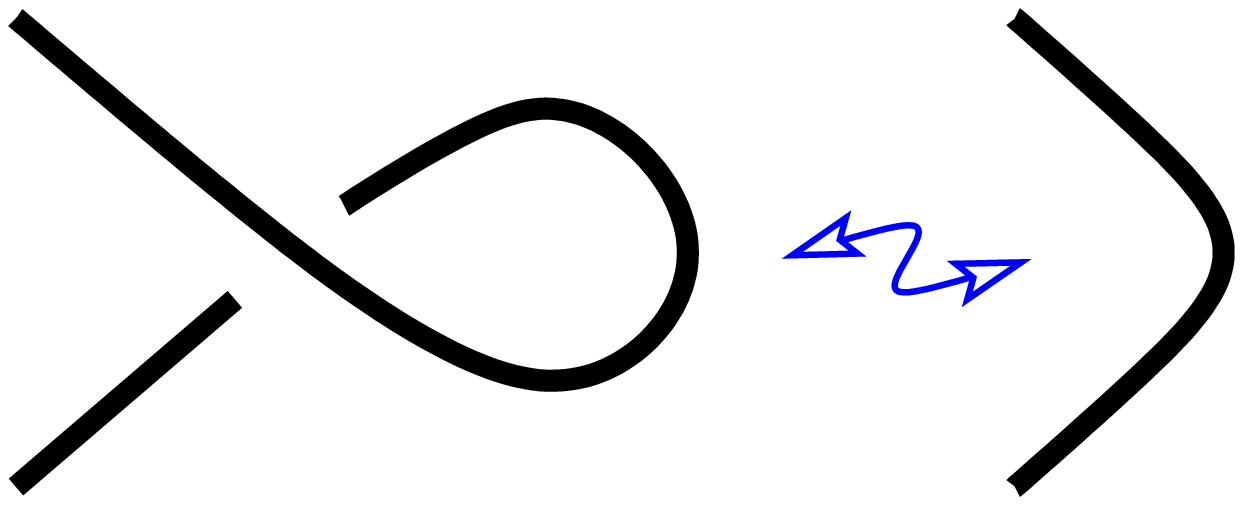}{}{75}{20}{20}\qquad\qquad
  \risS{-18}{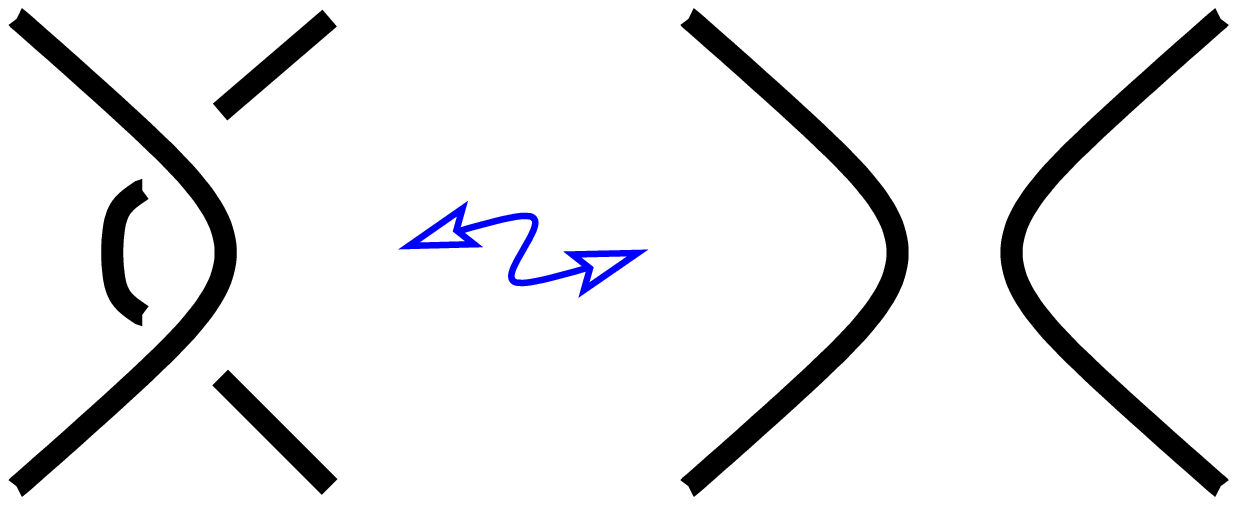}{}{75}{0}{20}\qquad\qquad
  \risS{-18}{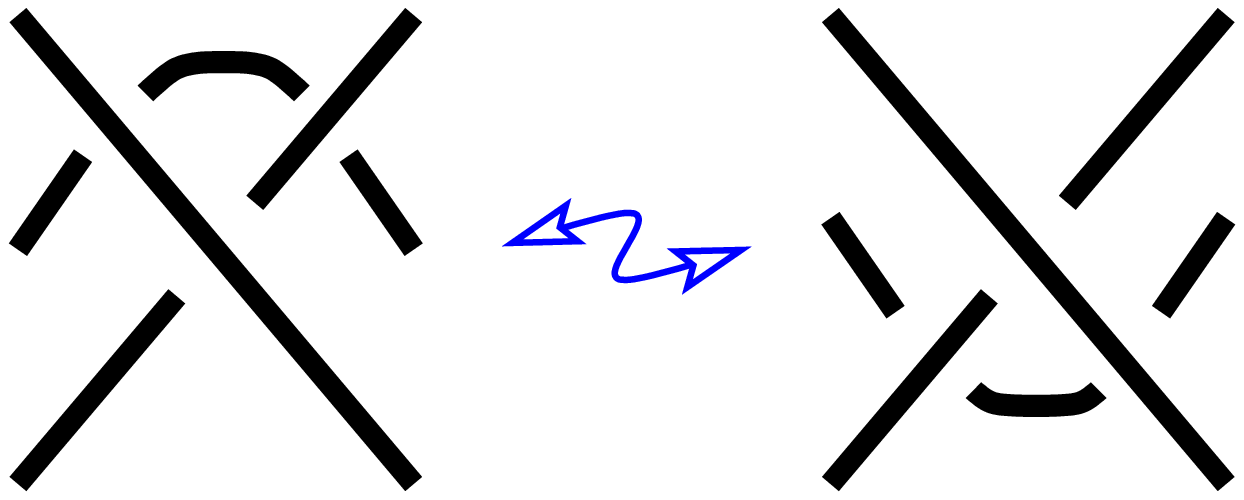}{}{75}{0}{20}
$$
and the {\it virtual} Reidemeister moves
$$\risS{-18}{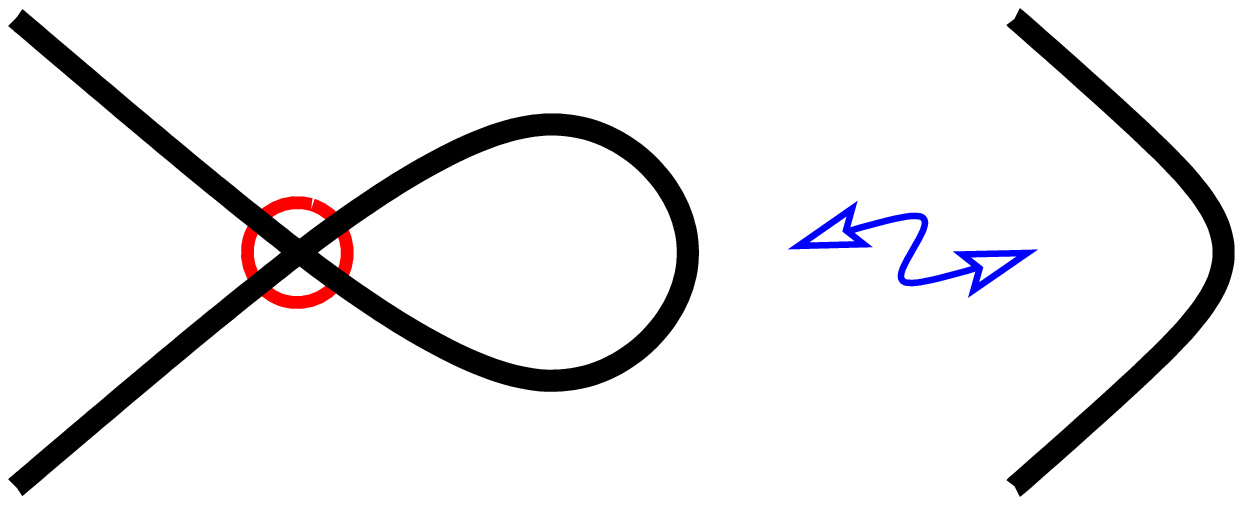}{}{75}{20}{20}\qquad\quad
  \risS{-18}{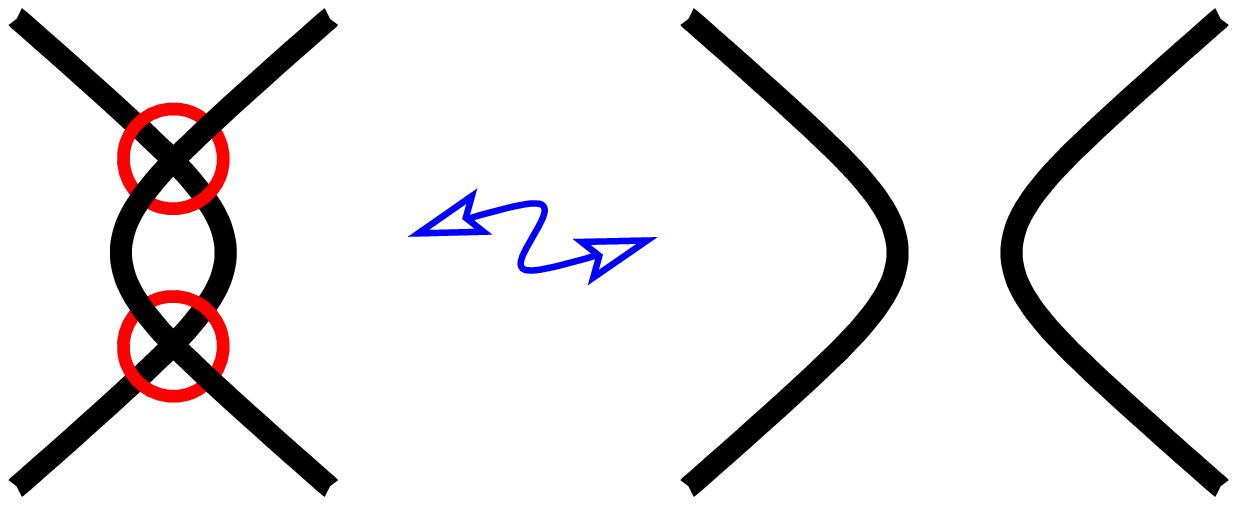}{}{75}{0}{20}\qquad\quad
  \risS{-18}{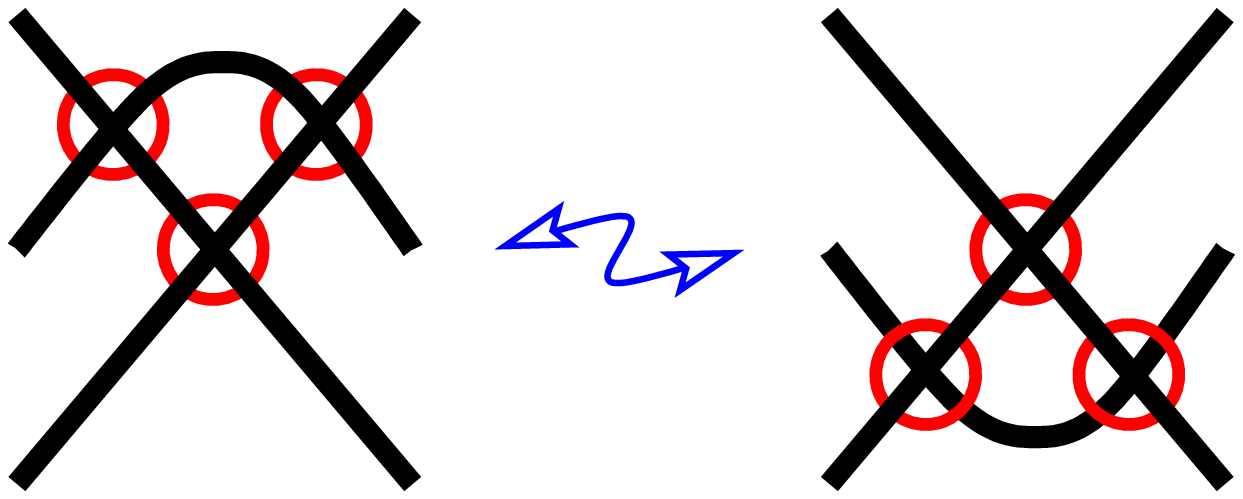}{}{75}{0}{20}\qquad\quad
  \risS{-18}{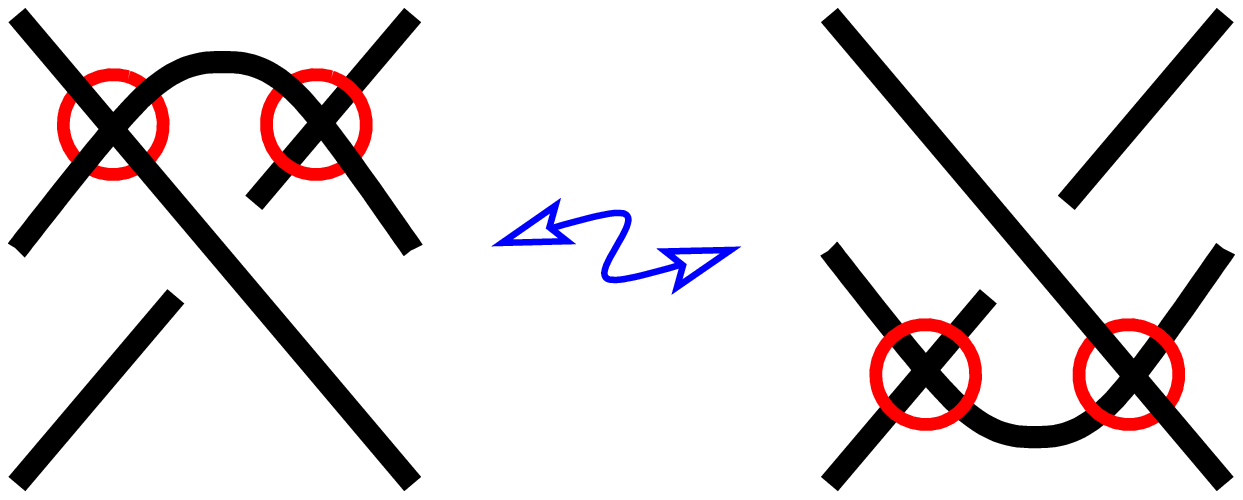}{}{75}{0}{20}
$$
Here the virtual crossings are encircled for the emphasis.

N.~Kamada introduced \cite{Kam1, Kam2} the notion of a {\it checkerboard coloring} of a
virtual link diagram. This is a coloring of one side of the diagram in its small
neighborhood, such that near a classical crossing it alternates like on a checkerboard,
and near a virtual crossing the colorings go through without noticing the crossing strand
and its coloring. Not every virtual link is checkerboard colorable. Here are two examples.
$$\risS{-40}{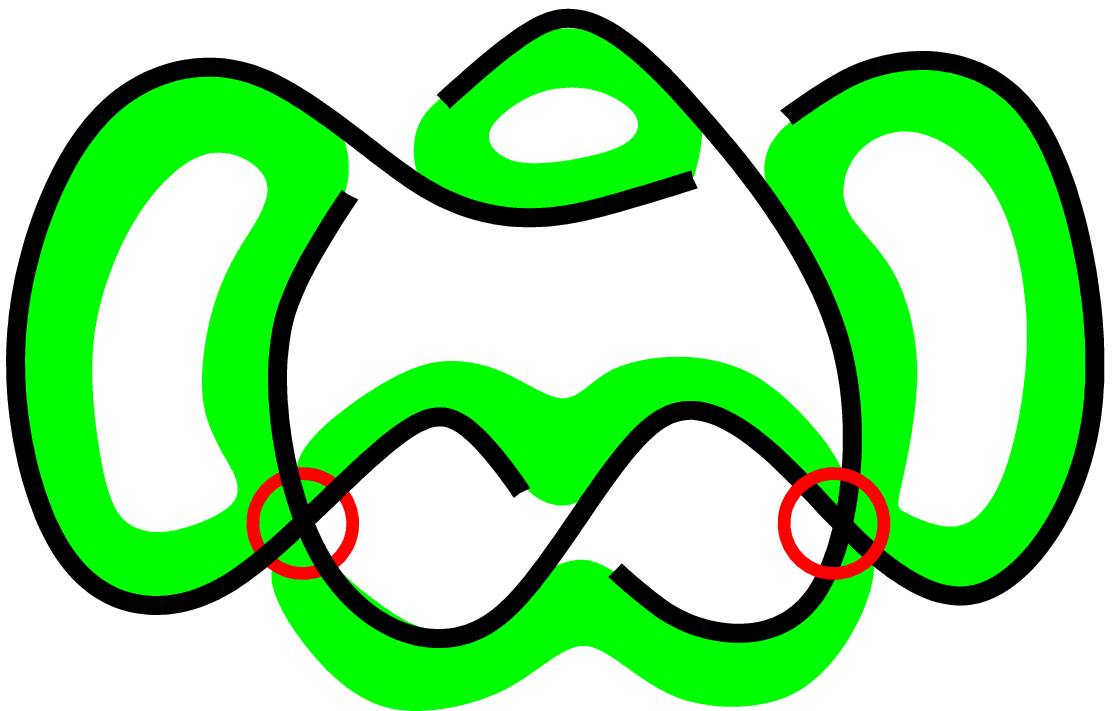}{\put(-8,-15){\mbox{\small checkerboard colorable}}}{80}{0}{65}
\hspace{150pt}\label{p-ex}
\risS{-35}{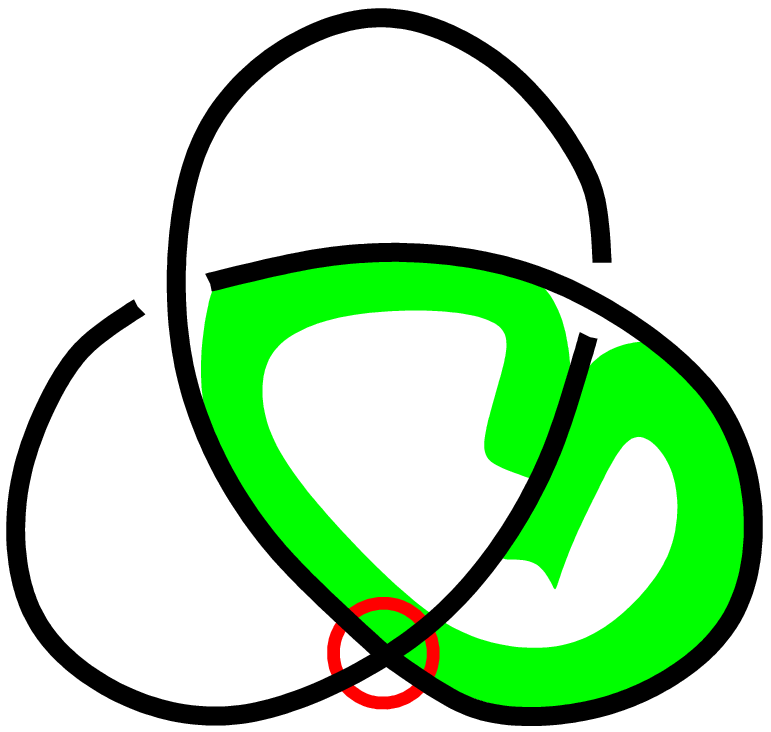}{\put(-40,-20){\mbox{\small not checkerboard colorable}}}{50}{0}{50}
$$
A similar notion was introduced and explored by V.~Manturov (see \cite{Man} and the
references therein), who called them {\it atoms}.

Note that the left virtual knot diagram above is alternating in a sense that classical
overcrossings and undercrossings alternate, while the right virtual diagram is not
alternating. It was proved in \cite{Kam1} that a virtual link diagram is checkerboard
colorable if and only if it can be made alternating by a number of classical crossing
changes. In particular, every classical link diagram is checkerboard colorable.

In \cite{Kam1, Kam2} N.~Kamada showed that many classical results on knots and links
can be extended to checkerboard colorable virtual links.  This paper is devoted to
a new result in this direction.
Namely, we generalize the celebrated theorem of M.~Thistlethwaite~\cite{Th}
(see also \cite{Ka1, Ka2}), which established a connection between the \emph{Jones
polynomial} for links and knots and the \emph{Tutte polynomial} for graphs.  Formally,
Thistlethwaite showed that, up to a sign and a power of~$t$, the Jones
polynomial $V_L(t)$ of an alternating link~$L$ is equal to the specialization of the
Tutte polynomial $T_{\Ga_L}(-t,-t^{-1})$ of the corresponding graph~$\Ga_L$. For
virtual links, the graph~$\Ga_L$ is naturally embedded into a surface rather than
into the plane, i.e. it becomes a \emph{ribbon graph}.
In this case, instead of the Tutte polynomial we should consider its generalization,
the \emph{Bollob\'as-Riordan polynomial}. Interestingly, the Bollob\'as-Riordan
polynomial was introduced with (very different) knot theoretic applications in
mind~\cite{BR2,BR3}.

The paper is structured as follows.  In the first two
sections we recall definitions of the \emph{Kauffman bracket} of virtual links
and the Bollob\'as-Riordan polynomial of ribbon graphs.
In section~3 we construct a ribbon graph from a checkerboard colorable
virtual link diagram and
state the Main Theorem for alternating virtual links.  As often appears in these cases,
the proofs of the results about (generalizations of) the Tutte
polynomial are quite straightforward. The proof of the
Main Theorem is postponed until section~5.  In section~4
we extend our results to signed ribbon graphs and derive
the Jones polynomial of an arbitrary checkerboard colorable
virtual link as an appropriate evaluation.  We conclude with final remarks
and an overview of the literature.

\bigskip
\section{The Kauffman bracket of virtual links.}
Let $L$ be a virtual link diagram.
Consider two ways of resolving a classical crossing.
The {\it $A$-splitting}, $\smf{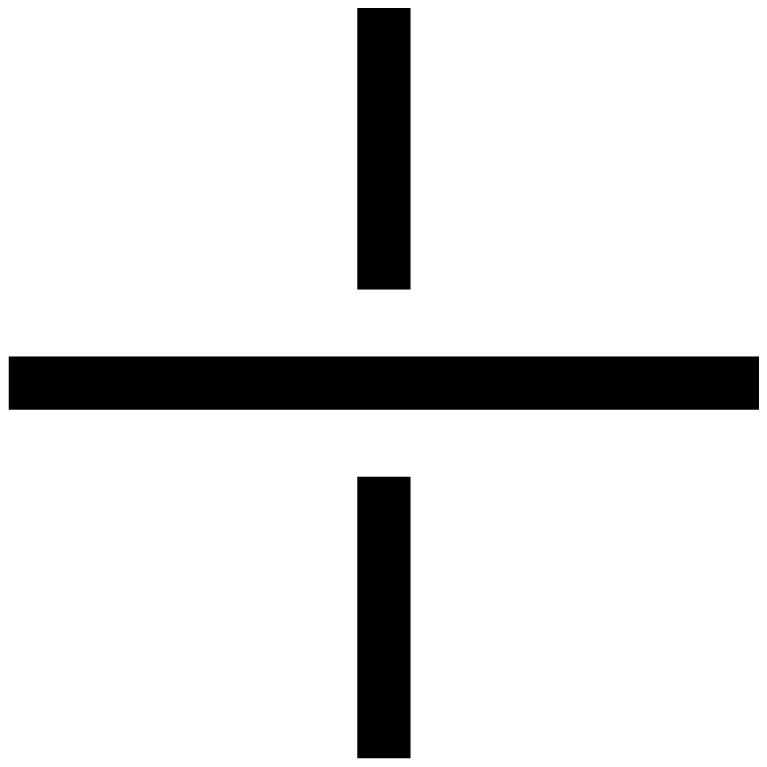}\ \leadsto\ \smf{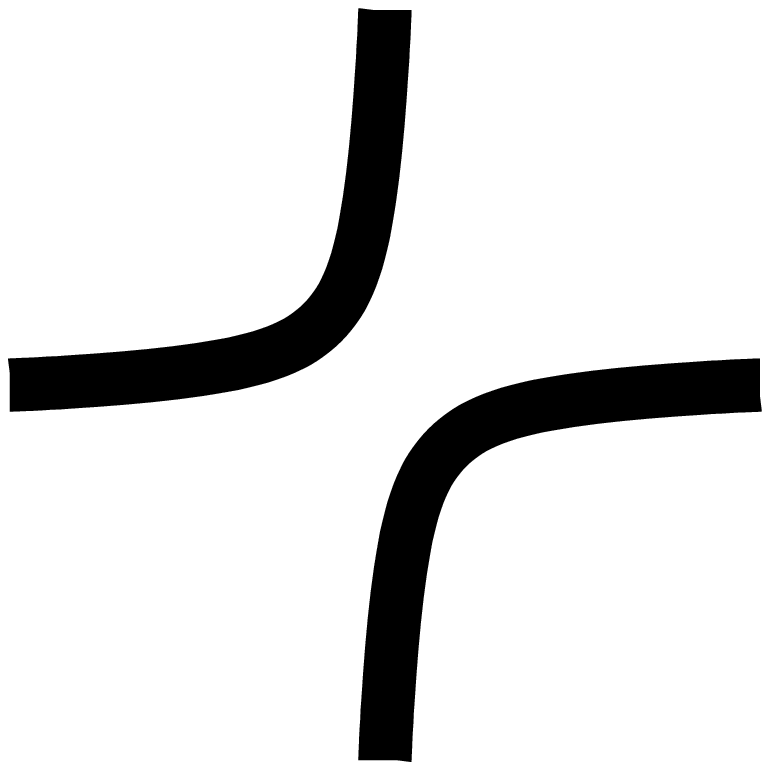}$,
is obtained by uniting the two regions swept out by the overcrossing arc under
the   counterclockwise rotation until the undercrossing arc.
Similarly, the {\it $B$-splitting},
$\smf{cr}\ \leadsto\ \smf{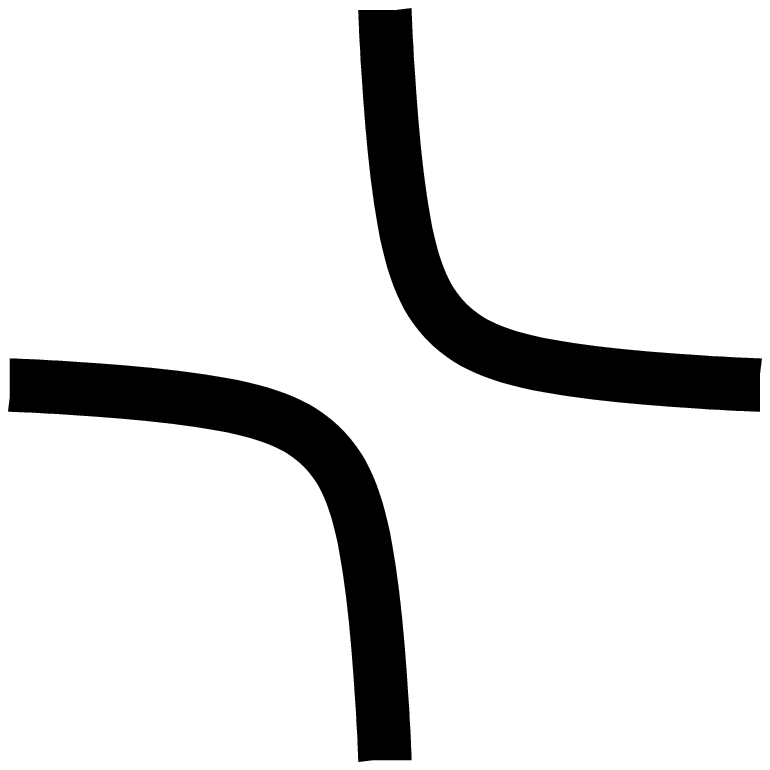}$, is
obtained by uniting the other two regions. A {\it state} $S$ of
a link diagram~$L$
is a way of resolving each classical crossing of the diagram.
Denote by $\cS(L)$ the set of the states of~$L$.
Clearly, a diagram $L$ with $n$ crossings has $|\cS(L)| = 2^n$
different states.

Denote by $\a(S)$ and $\b(S)$ the number of $A$-splittings and $B$-splittings
in a state $S$, respectively.  Also, denote by $\de(S)$ the number of
components of the curve obtained from the link
diagram $L$ by all
splittings according to the state $S \in \cS(L)$. Note that virtual crossings do not connect the components.

\begin{defn}\label{def:kb}
The \emph{Kauffman bracket} of a diagram $L$ is a polynomial in three variables
$A$, $B$, $d$ defined by the formula:
\begin{equation}
\kb{L} (A,B,d)\ :=\ \sum_{S \in \cS(L)} \,
A^{\a(S)} \, B^{\b(S)} \, d^{\de(S)-1}\,.
\end{equation}
\end{defn}

Note that $\kb{L}$ is \emph{not} a topological
invariant of the link and in fact depends on the link
diagram. However, it defines the \emph{Jones polynomial}
$J_L(t)$ by a simple substitution
$A=t^{-1/4}$, $B=t^{1/4}$, $d=-t^{1/2}-t^{-1/2}$:
$$J_L(t)\, := (-1)^{w(L)} t^{3w(L)/4} \kb{L} (t^{-1/4}, t^{1/4}, -t^{1/2}-t^{-1/2})\ .
$$
Here $w(L)$ denotes the {\it writhe}, determined
by the orientation of $L$ as the sum over the classical
crossings of $L$ of the following signs\,:
$$\begin{picture}(50,30)(0,0)
  \put(0,0){\includegraphics[width=50pt]{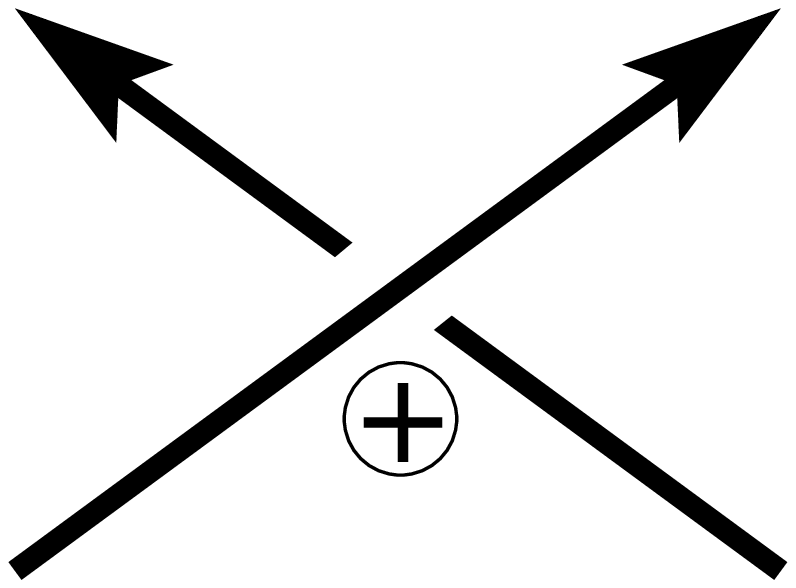}}
     \end{picture} \hspace{3cm}
\begin{picture}(50,40)(0,0)
  \put(0,0){\includegraphics[width=50pt]{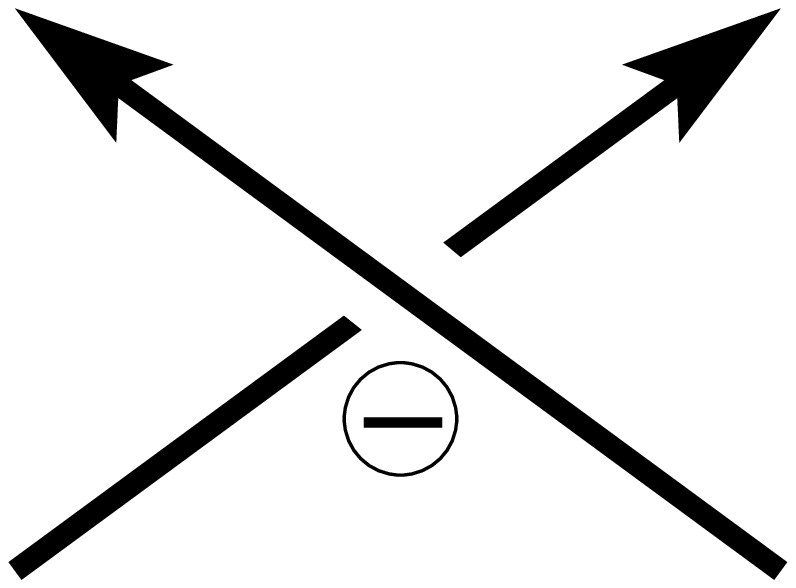}}
     \end{picture}
$$
The Jones polynomial is a classical topological invariant (see e.g.~\cite{B,Man}).

\begin{Example}\label{ex1} {\rm
Consider the virtual knot diagram
$L$ from the example above and shown on the left of the table below.
It has two virtual and three classical crossings, so there are
eight states for it,
$|\cS(L)| = 8$.
The curves obtained by the splittings and the corresponding
parameters $\a(S)$, $\b(S)$, and $\de(S)$ are shown in the
remaining columns of the table.
$$
\begin{array}{c||c|c|c|c}
\laf{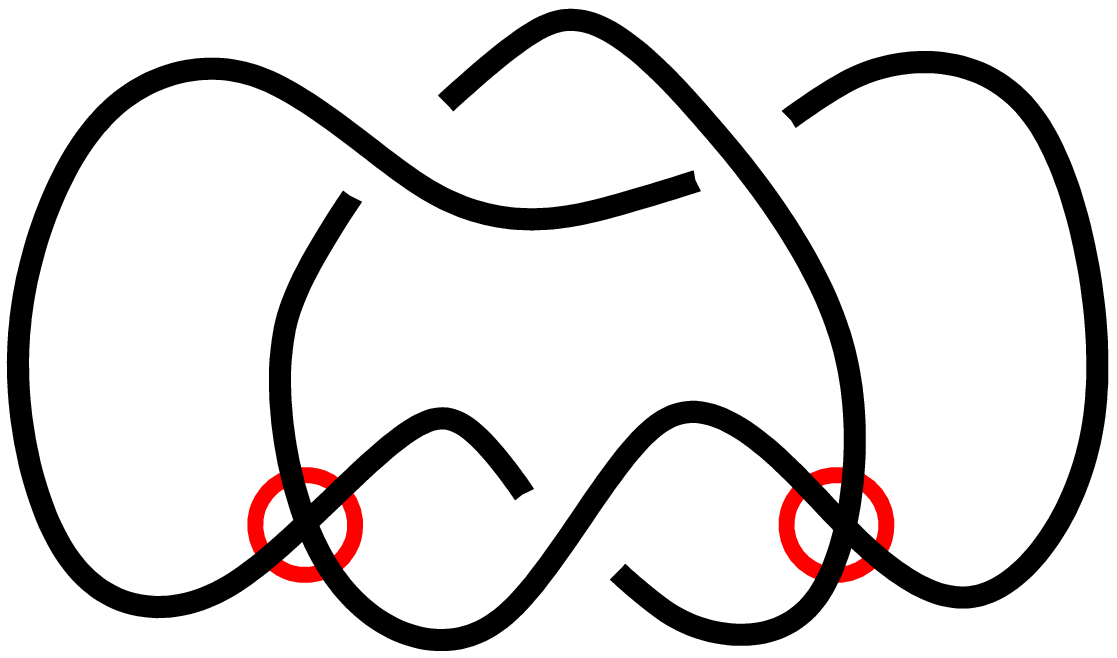} & \laf{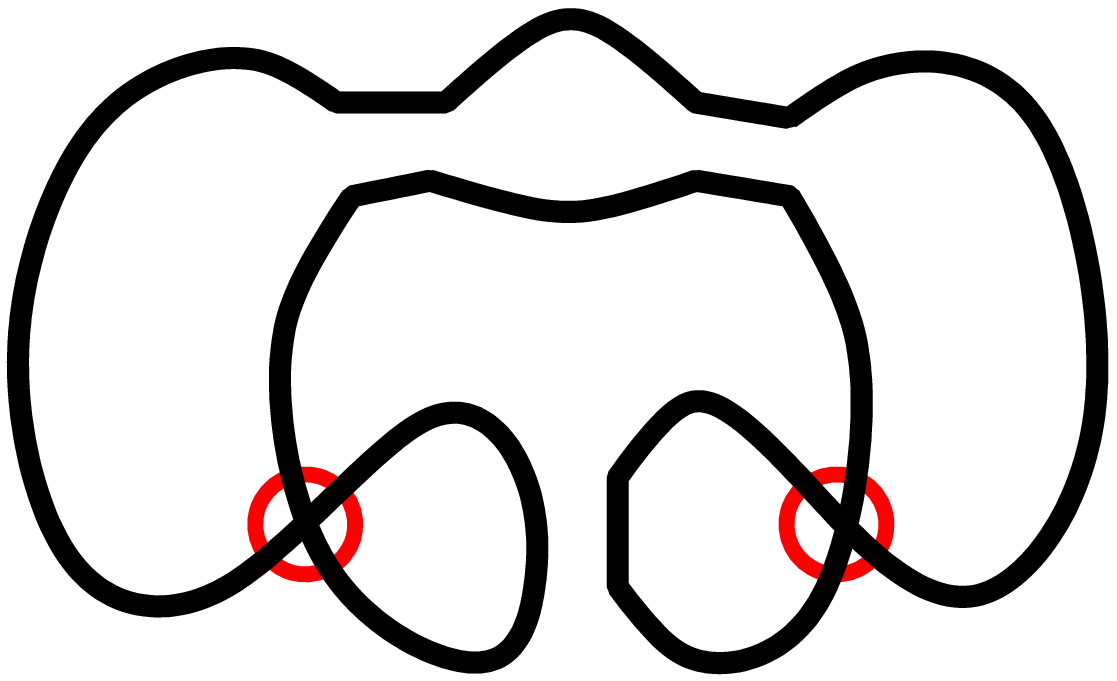} & \laf{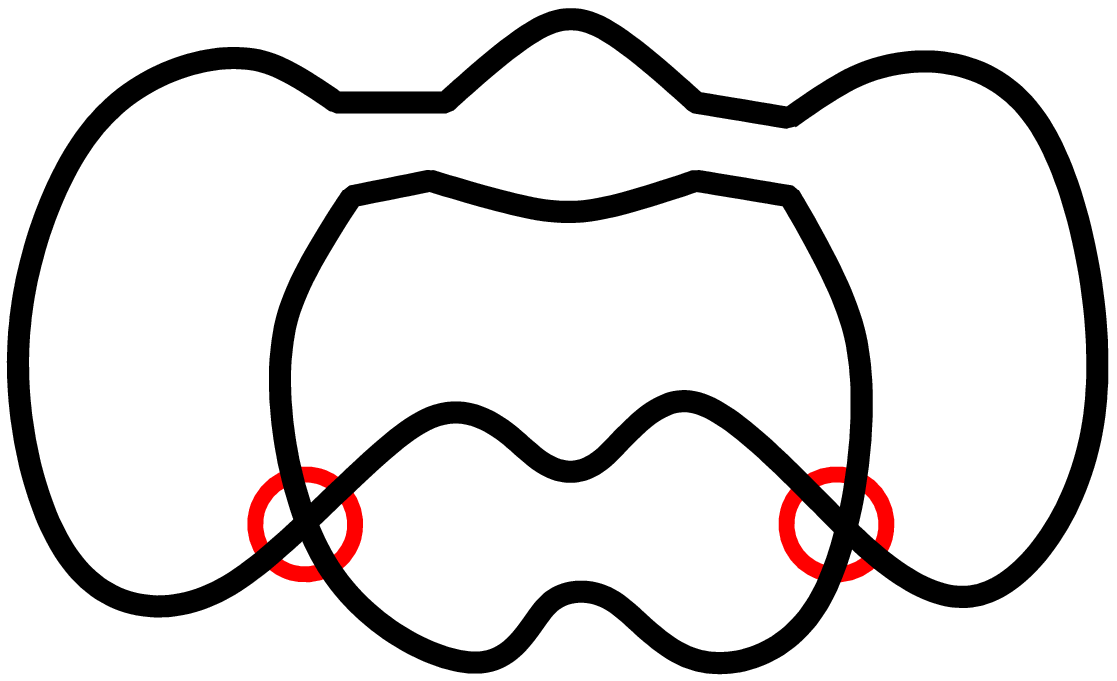} & \laf{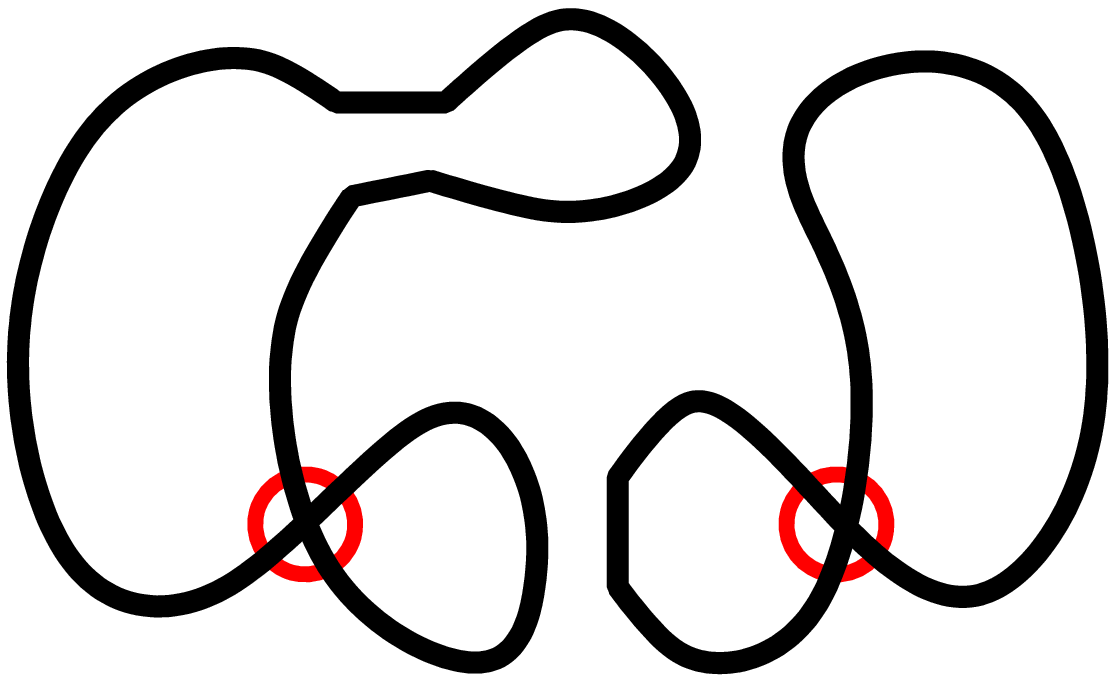} & \laf{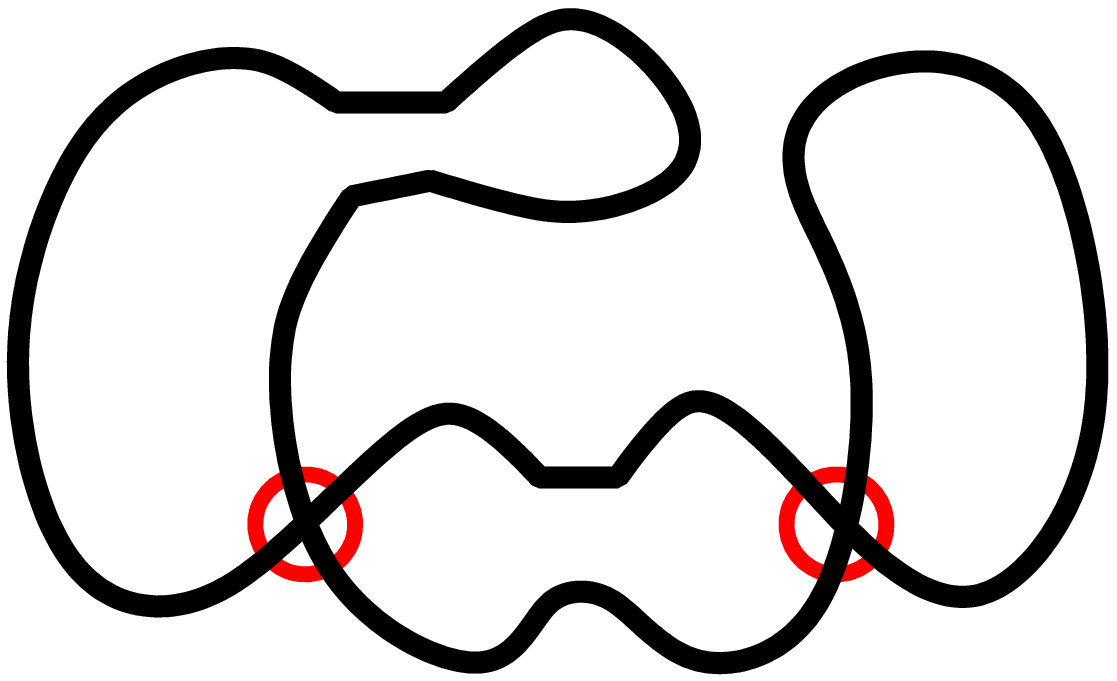}\\ \hline
 (\a,\b,\de) & (3,0,1) & (2,1,2) & (2,1,2) & (1,2,1)\makebox(0,15){}
\\ \hline\hline
& \laf{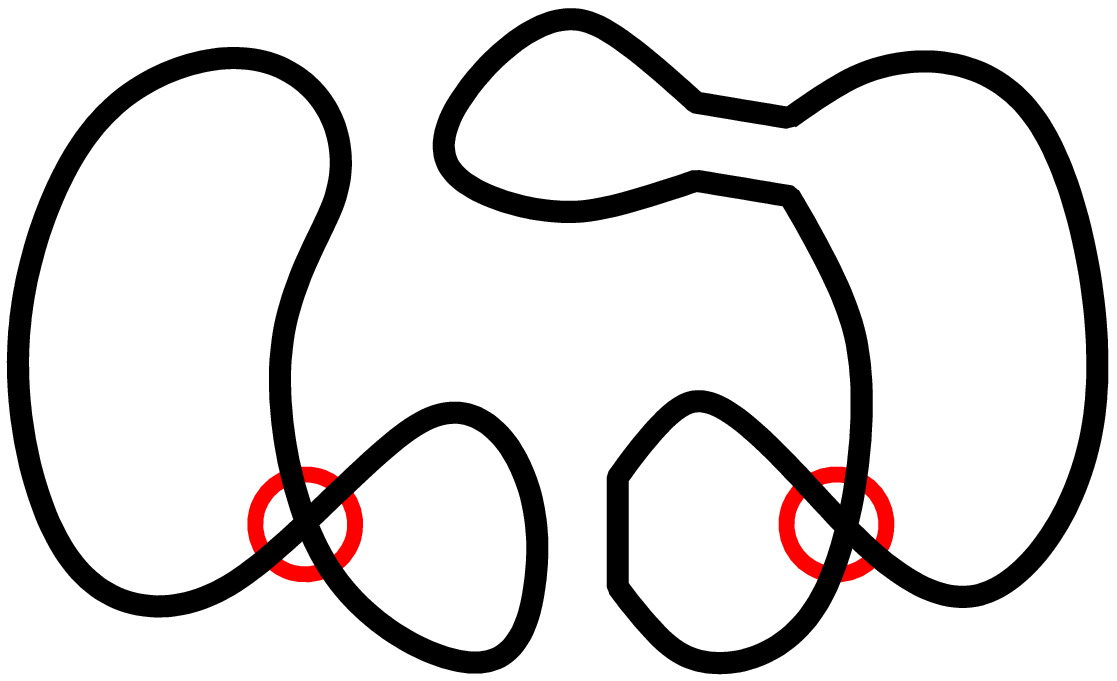} & \laf{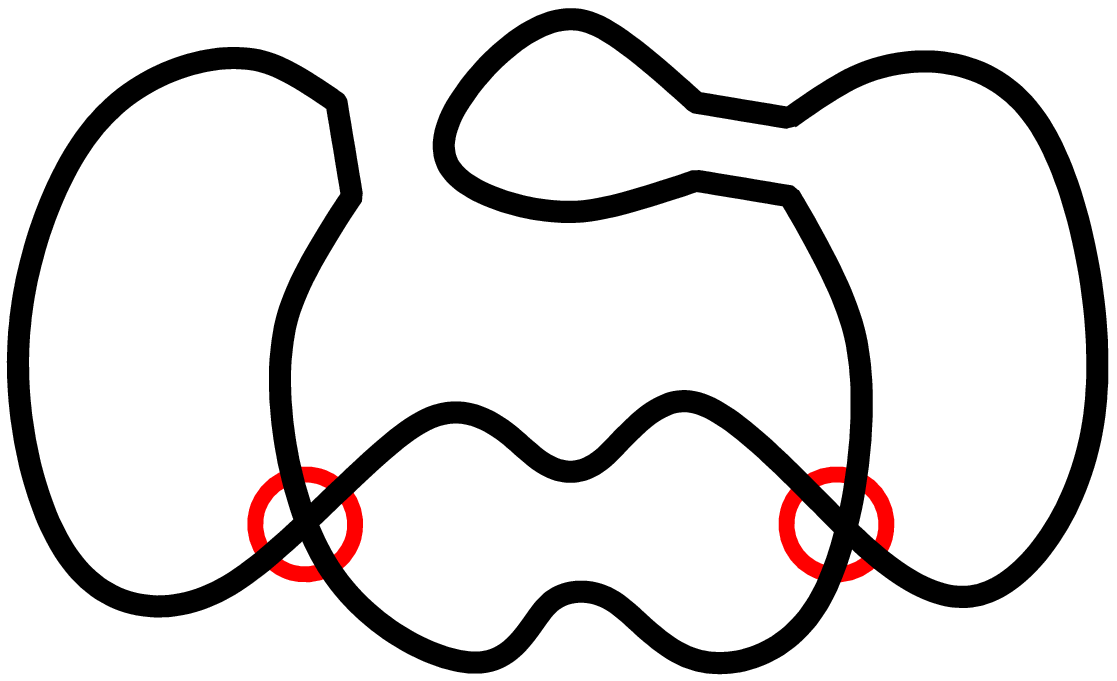} & \laf{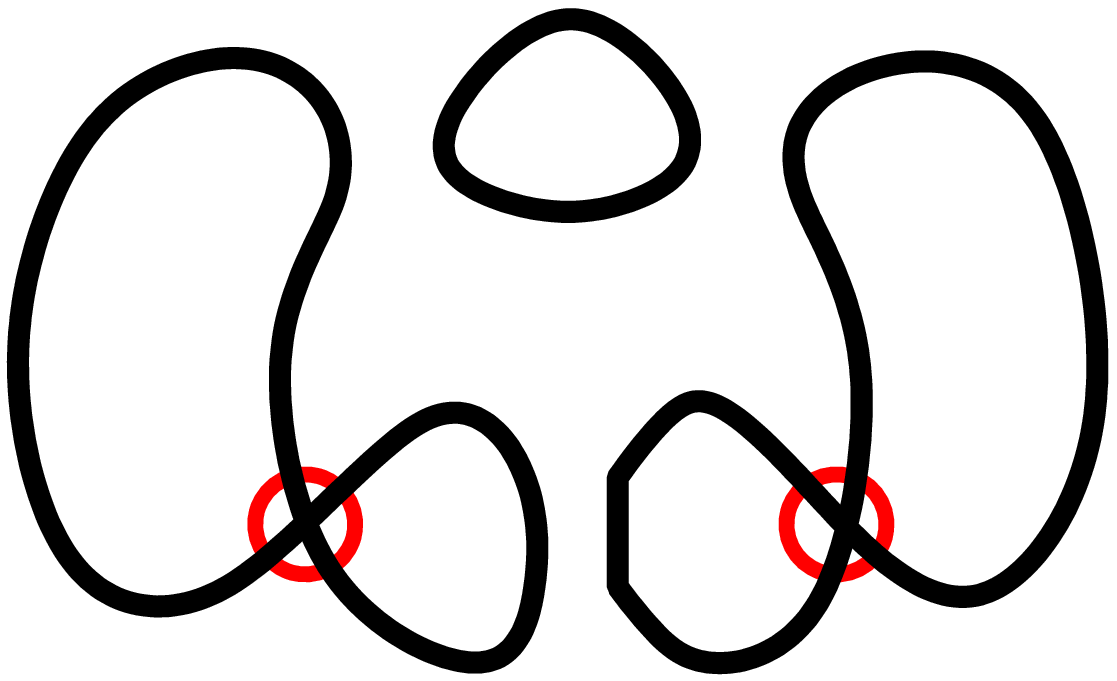} & \laf{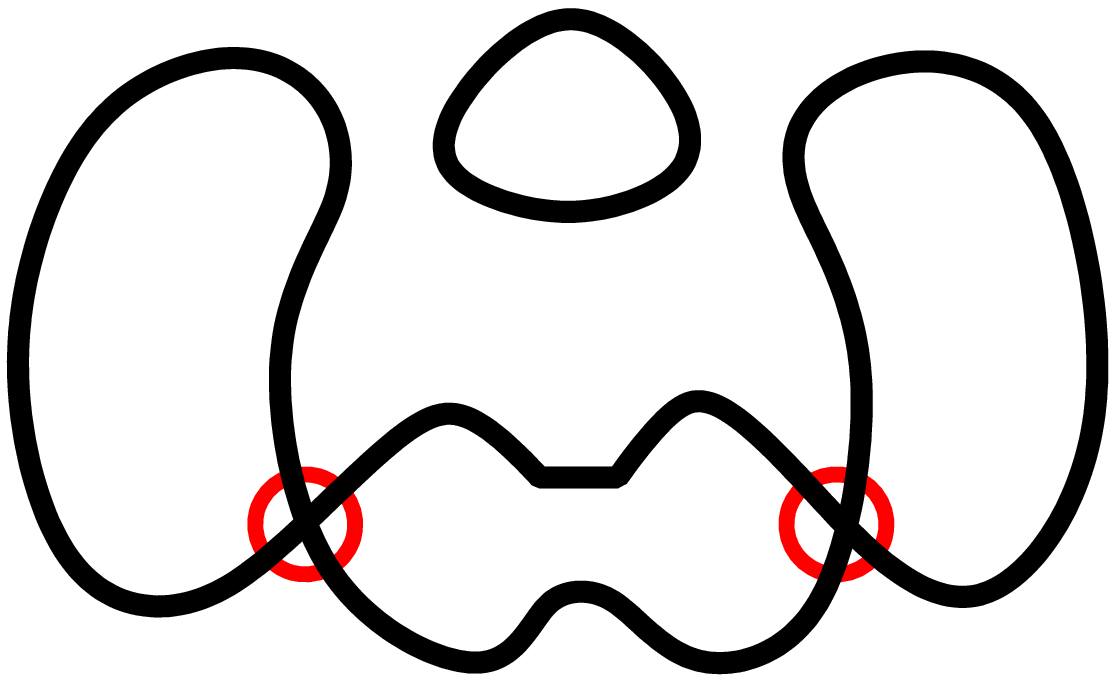}\\ \cline{2-5}
& (2,1,2) & (1,2,1) & (1,2,3) & (0,3,2)\makebox(0,15){}
\end{array}
$$

In this case the Kauffman bracket of $L$ is given by
$$\kb{L} = A^3 + 3A^2Bd + 2AB^2 + AB^2d^2+ B^3d\ .$$
It is easy to check that $J_L(t)=1$.
}
\end{Example}

\bigskip

\section{The Bollob\'as-Riordan polynomial.}

Let $\Ga = (V,E)$ be a undirected graph with the set of vertices $V$
and the set of edges $E$ (loops and multiple edges are allowed).
Suppose in each vertex $v \in V$ there is a fixed cyclic order on
edges adjacent to $v$ (loops are counted twice).  We call this
combinatorial structure a {\it ribbon graph}, and denote it
by $G$.  One can represent $G$ by making vertices into oriented `discs'
and connecting them by oriented `ribbons' as prescribed by the cyclic
orders (see Example~\ref{ex2} below) and in compliance with the orientation.
This defines a 2-dimensional surface with boundary, which by a
slight abuse of notation we also denote by $G$.

Formally, $G$ is the surface with boundary represented as
the union of two sets of closed topological discs, corresponding to
vertices $v \in V$ and edges $e \in E$, satisfying the following
conditions:

\smallskip

$\bu$ \ these discs and ribbons intersect by disjoint line segments,

$\bu$ \ each such line segment lies on the boundary of precisely
\\ \hspace*{.75cm} one vertex and precisely one edge,

$\bu$ \  every edge contains exactly two such line segments.

\smallskip

\noindent
It will be clear from the context whether by $G$ we mean the
ribbon graph or its underlying surface.  In this paper we
restrict ourselves to oriented surfaces~$G$.
We refer to \cite{GT} for other definitions and references.

For a ribbon graph $G$, let
$v(G) = |V|$ denote the number of vertices,
$e(G) = |E|$ denote the number of edges, and
$k(G)$ denote the number of connected components of~$G$.
Also, let $r(G)=v(G)-k(G)$ be the {\it rank} of~$G$, and
$n(G)=e(G)-r(G)$ be the {\it nullity} of~$G$.
Finally, let $\bc(G)$ be the number of connected
components of the boundary of the surface~$G$.

A spanning subgraph of a ribbon graph $G$ is defined as a
subgraph which contains all the vertices, and a subset of the edges.
Let $\cF(G)$ denote the set of the spanning subgraphs of~$G$.
Clearly, $|\cF(G)| = 2^{e(G)}$.

\begin{defn} \label{def:br}
The \emph{Bollob\'as-Riordan polynomial} $R_G(x,y,z)$ of a ribbon
graph $G$ is defined by the formula
\begin{equation}\label{def_br}
R_G(x,y,z)\ :=\ \sum_{F \in \cF(G)}
   x^{r(G)-r(F)} y^{n(F)} z^{k(F)-\bc(F)+n(F)}\,.
\end{equation}
\end{defn}

This version of the polynomial is obtained from the original
one \cite{BR2,BR3} by a simple substitution $x+1$ for $x$.
Note that for a planar ribbon graph~$G$
(i.e. when the surface~$G$ has genus zero)
the Euler's formula gives $k(F)-\bc(F)+n(F)=0$ for all $F\subseteq G$.
Therefore, the Bollob\'as-Riordan polynomial~$R_G$ does not
contain powers of~$z$. In fact, in this case it is essentially
equal to the classical Tutte polynomial $T_{\Ga}(x,y)$ of the
(abstract) core graph~$\Ga$ of~$G$:
$$R_G(x-1,y-1,z) = T_{\Ga}(x,y)\,.$$
Similarly, a specialization $z=1$ of the Bollob\'as-Riordan polynomial
of an arbitrary ribbon graph~$G$, gives the Tutte polynomial once
again:
$$R_G(x-1,y-1,1) = T_{\Ga}(x,y)\,.$$
We refer to~\cite{BR2,BR3} for proofs of these formulas and
to~\cite{B,W} for general background on the Tutte polynomial.

\begin{Example}\label{ex2} {\rm
Consider the ribbon graph $G$ shown on the left in the table below.
The other columns show eight possible spanning subgraphs~$F$ and
the corresponding values of $k(F)$, $r(F)$, $n(F)$ and $\bc(F)$.
$$
\begin{array}{c||c|c|c|c}
\baf{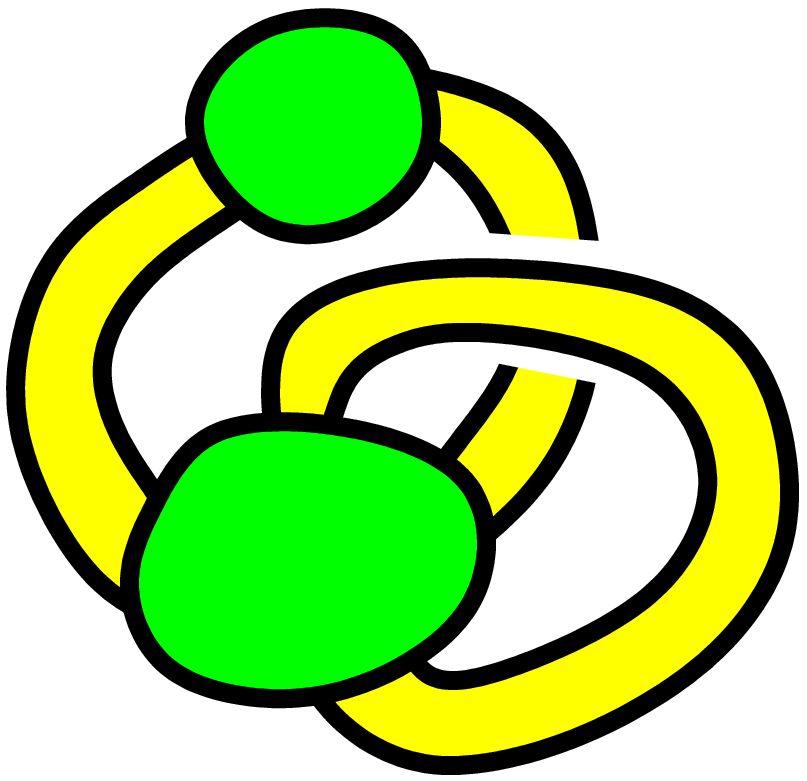} & \baf{ex2AAA} & \saf{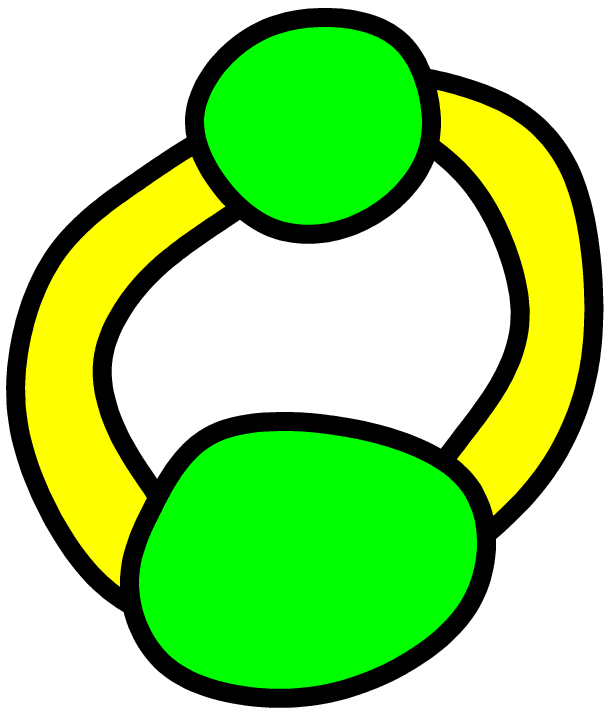}{40} & \baf{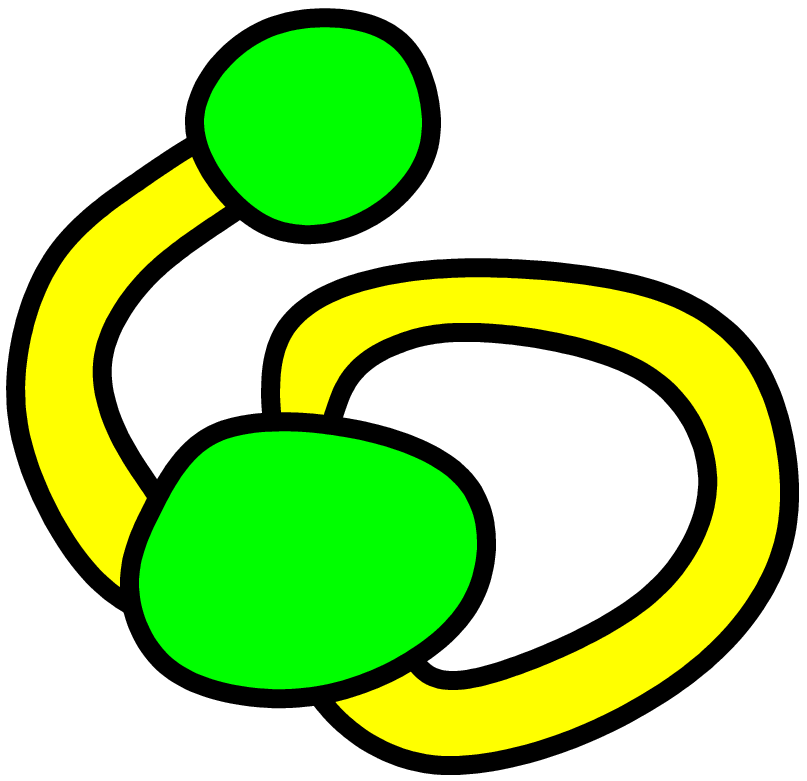} & \saf{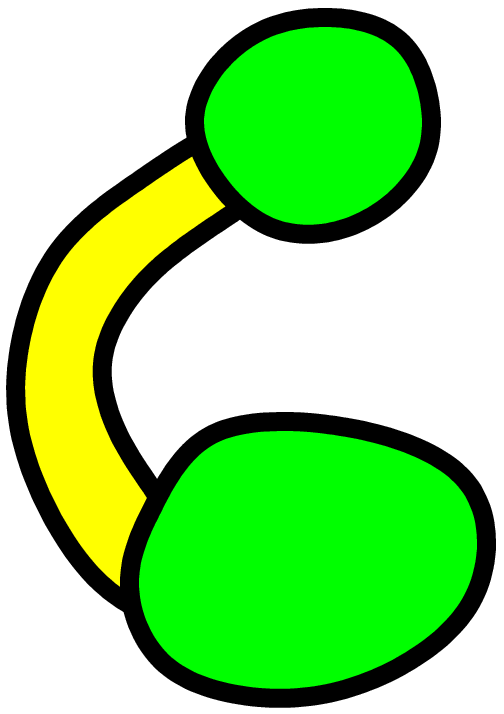}{32}\\ \hline
 (k,r,n,\bc) & (1,1,2,1) & (1,1,1,2) & (1,1,1,2) & (1,1,0,1)\makebox(0,15){}
\\ \hline\hline
& \rb{-5}{\saf{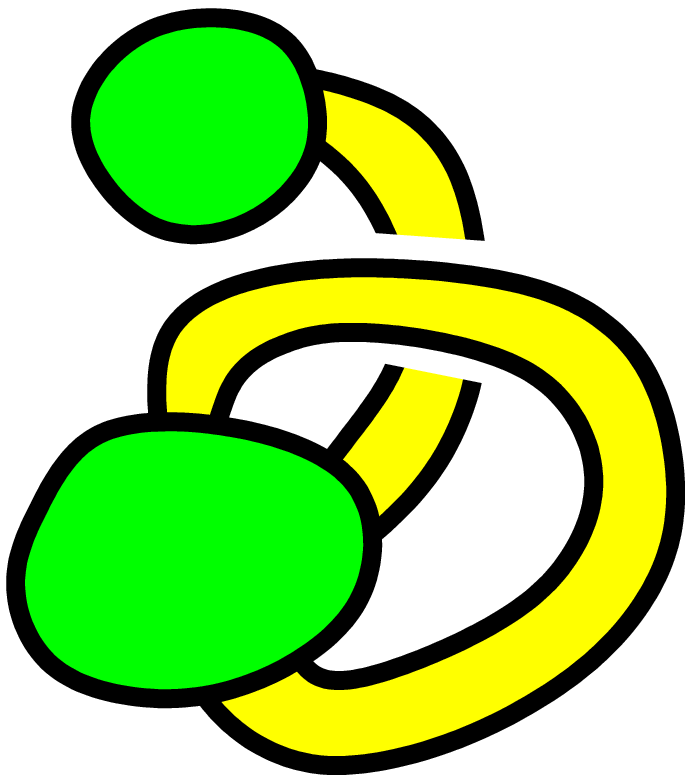}{45}} & \saf{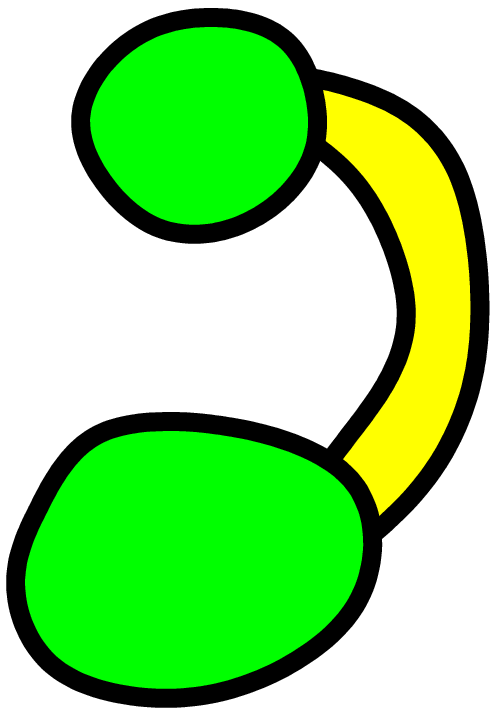}{32} & \
  \rb{-5}{\saf{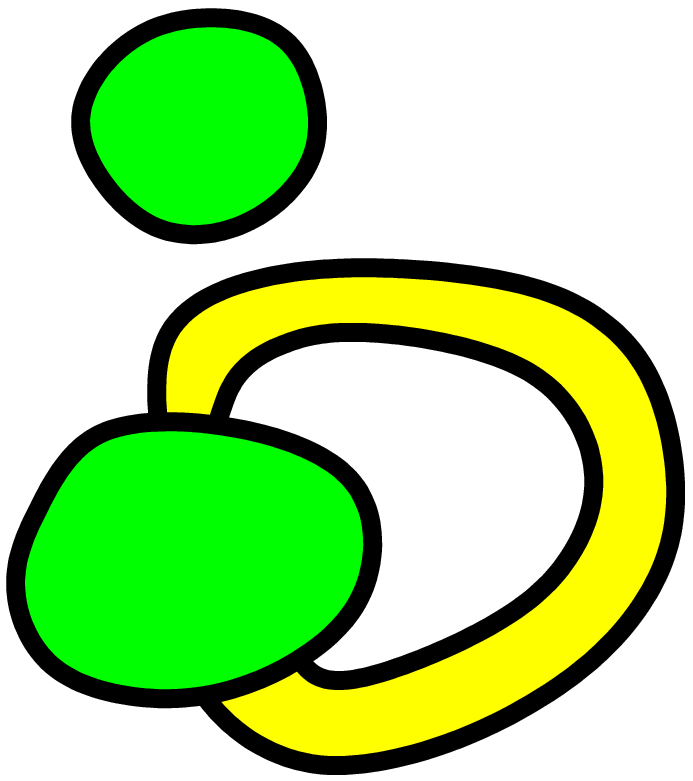}{45}} & \quad\saf{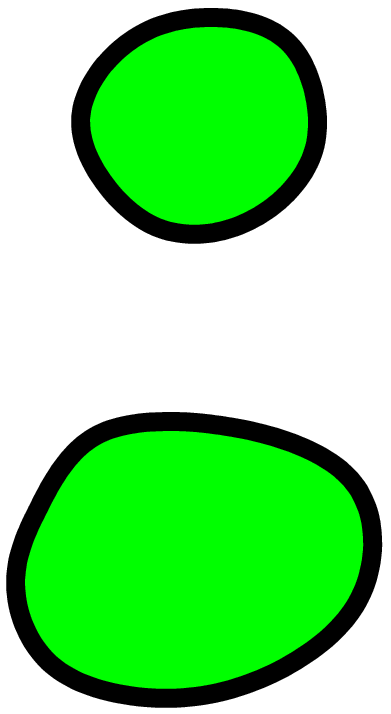}{25}\\ \cline{2-5}
& (1,1,1,2) & (1,1,0,1) & (2,0,1,3) & (2,0,0,2)\makebox(0,15){}
\end{array}
$$

\bigskip
\noindent
Now use the definition to compute the corresponding
Bollob\'as-Riordan polynomial:
$$R_G(x,y,z) \, = \, y^2z^2 + 3y + 2 + xy + x$$
}
\end{Example}

\bigskip
\section{Ribbon graphs from virtual diagrams and the Main Theorem.}

In this section we construct a ribbon graph $G_L$ starting with an
alternating virtual link diagram~$L$.

\smallskip
It was shown in \cite{Kam2} that every alternating link diagram~$L$ has
a canonical checkerboard coloring which can be constructed in the following way.
Near every classical crossing we color black the vertical angles swept out
by the overcrossing arc under the counterclockwise rotation, i.e. angles
that are glued together by the $A$-splitting.  Since the diagram is
alternating, it is evident that all these local colorings near classical
crossings agree with each other and can be extended to a global checkerboard
coloring. For example, the checkerboard coloring of the virtual knot
in the introduction is canonical.

Topologically a coloring is represented by a bunch of annuli. Each annulus has
two boundary circles, an exterior circle which goes along the link except small
arcs near classical crossings where it jumps from one strand to another one,
and an interior circle. In order to construct a ribbon graph from a (virtual)
link diagram we replace every crossing by an edge-ribbon connecting the
corresponding arcs of the exterior circles:
$$\risS{-20}{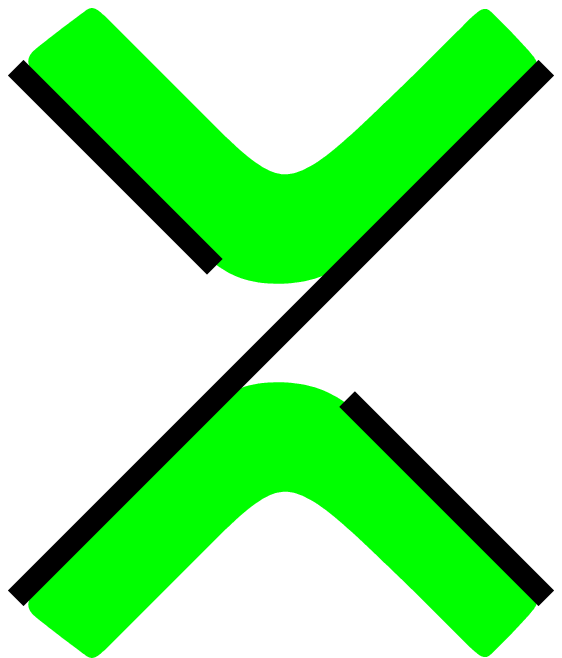}{}{40}{25}{20}
\quad\risS{-5}{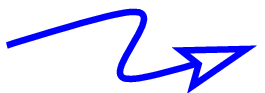}{}{25}{0}{20}\quad
\risS{-20}{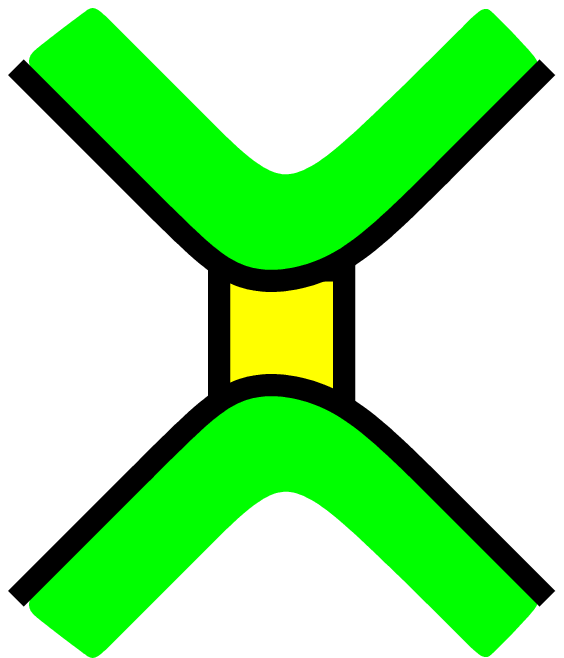}{}{40}{25}{20}
\hspace{100pt}
\risS{-20}{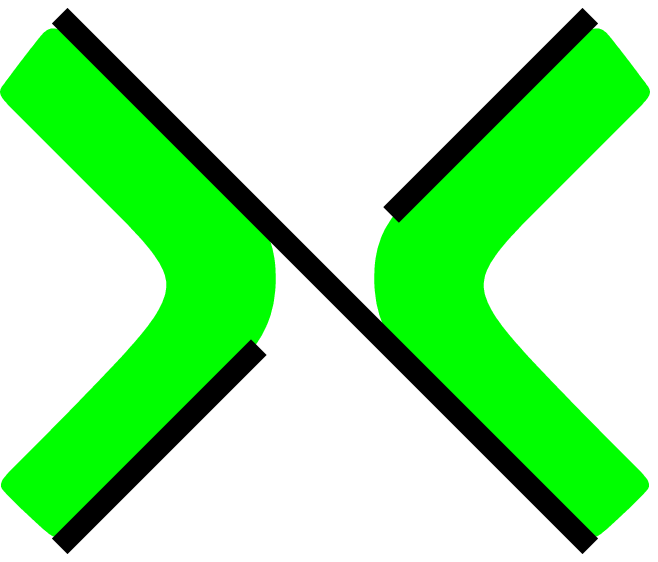}{}{50}{25}{20}
\quad\risS{-5}{to}{}{25}{0}{20}\quad
  \risS{-20}{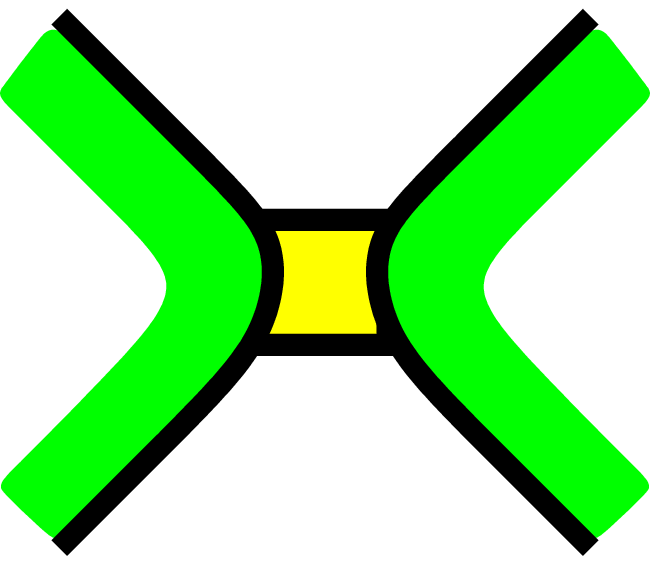}{}{50}{25}{20}
$$
The ribbon graph $G_L$ is obtained by gluing discs along the interior circles
of the annuli of the coloring. Here is a ribbon graph constructed
from a virtual link diagram in the introduction:
$$\risS{-25}{col-kn}{}{80}{20}{25}
\qquad\risS{-5}{to}{}{25}{0}{20}\qquad
\risS{-25}{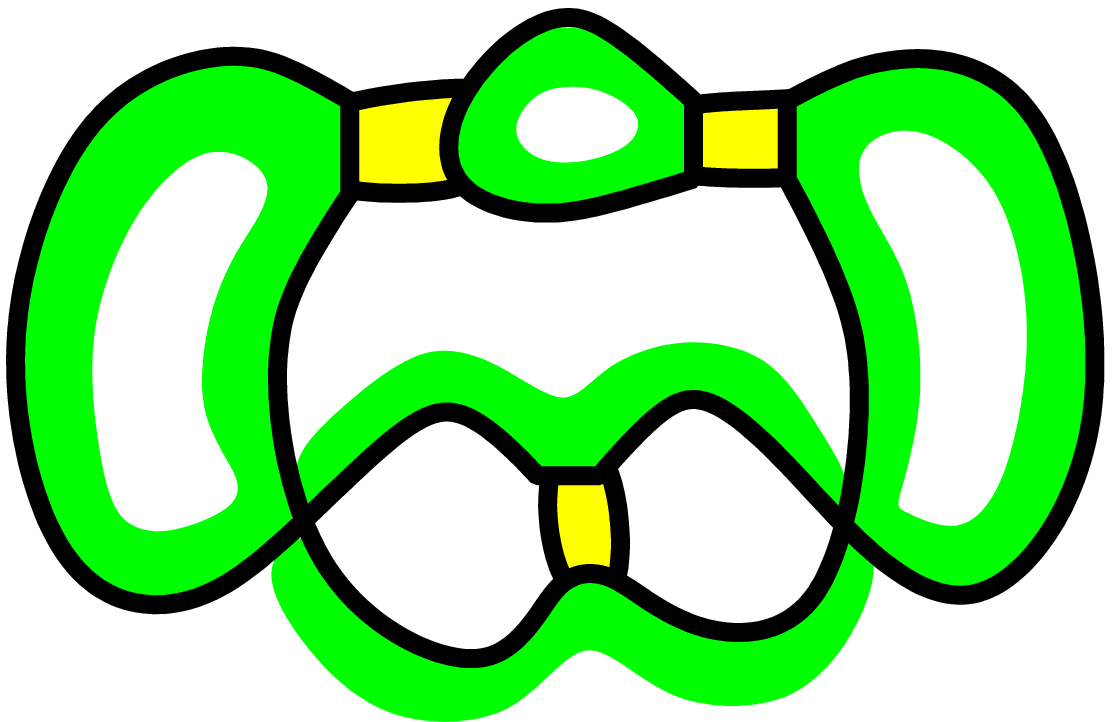}{}{80}{20}{20}\ = \
\risS{-35}{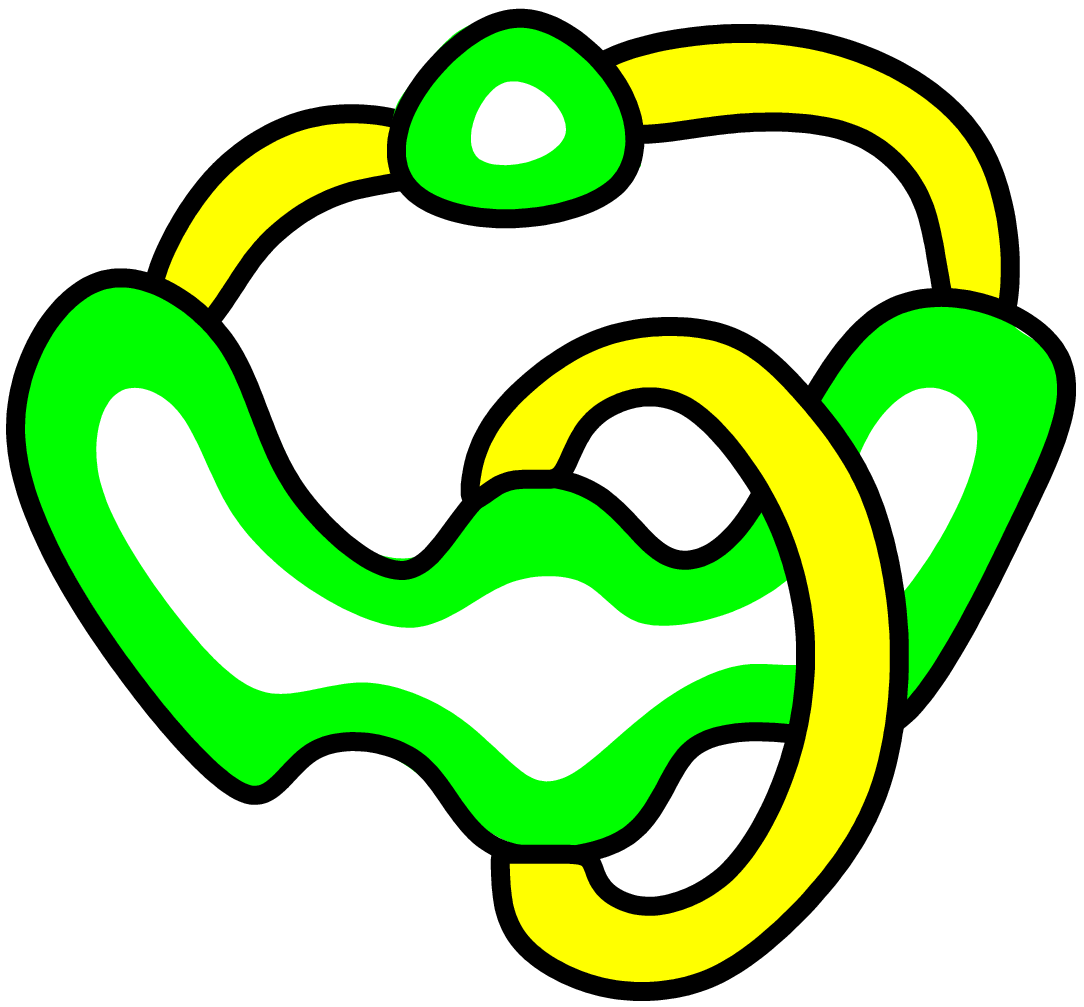}{}{80}{30}{40}\ = \
\risS{-20}{ex2AAA}{}{50}{0}{20}
$$
Note that here we ignore the way a ribbon graph is embedded into
$\mathbb{R}^3$
and consider it as an abstract 2-dimensional surface.  Of course, every ribbon graph
can be obtained in this way from an appropriate virtual link diagram.

\medskip

\noindent
{\bf Main Theorem 3.1.}  {\it Let $L$ be an alternating virtual link
diagram and $G_L$ be the corresponding ribbon graph. Then
$$\kb{L} (A,B,d) \, = \, A^{r(G)} B^{n(G)} d^{k(G)-1} \,
R_{G_L}\left(\frac{Bd}{A}, \frac{Ad}{B}, \frac{1}{d}\right).$$
}

\medskip

We should warn the reader that although both the Kauffman bracket
and the Bollob\'as-Riordan polynomial are polynomials in three
variables, the former has only two free variables since
$\kb{L}$ is always homogeneous in~$A$ and~$B$.
This follows from the identity $\al(S) + \be(S) = e(S)$ for
all $S \in \cS(L)$.  Another way to see this is to note that
the values at which the Bollob\'as-Riordan polynomial $R_G(x,y,z)$
is evaluated in the theorem satisfy the equation $xyz^2 = 1$.
Thus, the situation here is different from the planar case where the
Kauffman bracket $\kb{L_G}(A,B,d)$ and the Tutte polynomial
$T_{\Ga}(x,y)$ determine each other.

\bigskip
\section{Extensions and applications.}

Define a \emph{signed ribbon graph} $\wh G$ to be a ribbon
graph $G$ given by $(V,E)$, and a \emph{sign function}
$\ve: E \to \{\pm 1\}$.
For a spanning subgraph $F \subseteq \wh G$ denote by $e_{-}(F)$
the number of edges $e \in E$ with $\ve(e)=-1$.
Denote by $\wo F = \wh G - F$ the complement to $F$ in $\wh G$, i.e.
the spanning subgraph of $\wh G$ with only those (signed) edges of~$G$
that do not belong to $F$.  Finally, let
$$s(F)  = \frac{e_{-}(F)-e_{-}(\wo F)}{2}$$
We define the
\emph{signed Bollob\'as-Riordan polynomial} $R_{\wh G}(x,y,z)$
as follows:
\begin{equation}\label{def_sbr}
R_{\wh G}(x,y,z)\ :=\ \sum_{F \in \cF(\wh G)}
   x^{r(\wh G)-r(F)+s(F)}
   y^{n(F)-s(F)}
   z^{k(F)-\bc(F)+n(F)}\, .
\end{equation}
This is a Laurent polynomial in $x^{1/2}$, $y^{1/2}$, and $z$.

\smallskip
Any checkerboard colorable virtual link diagram $L$ can be made alternating $\wt L$ by switching some classical overcrossings to undercrossings \cite{Kam1}.
We can label the
edges of $G_{\wt L}$ corresponding to the crossings where the switching was performed by $-1$, and the other edges by $+1$. The result is a signed ribbon graph denoted by
${\wh G}_L$. Note that it is not uniquely determined by the diagram $L$ because we can switch different sets of crossings to get an alternating diagram. However, the next theorem does not depend on this ambiguity.

\medskip
\begin{thm} \label{t:signed}
Let $L$ be a checkerboard colorable virtual link diagram, and ${\wh G}_L$ be
the corresponding signed ribbon graph. Then
$$\kb{L} (A,B,d)  \ = \ A^{r({\wh G}_L)} B^{n({\wh G}_L)} d^{k({\wh G}_L)-1} \,
R_{{\wh G}_L}\left(\frac{Bd}{A}, \frac{Ad}{B}, \frac{1}{d}\right).$$
\end{thm}

\smallskip

The proof follows verbatim the proof of the Main Theorem
(see the next section).  We leave the details to the reader.

\smallskip
The following result is an immediate consequence of Theorem~\ref{t:signed}.

\medskip

\begin{cor} \label{c:jones}
Let $\wh G$ be a signed ribbon graph corresponding
to a  checkerboard colorable virtual link diagram $L$ endowed
with an orientation. Then
$$J_L(t) \, = \, (-1)^{w(L)} \,
t^{\frac{3w(L)-r(\wh G)+n(\wh G)}{4}} \,
      \bigl(-t^{1/2}\!-t^{-1/2}\bigr)^{k(\wh G)-1} \,
R_{\wh G}\Bigl(-t-1, -t^{-1}\!-1, \frac{1}{-t^{1/2}\!-t^{-1/2}} \Bigr).
$$
In particular, if $\wh G$ is a planar ribbon graph with only positive edges
and $\Ga$ is its core graph, we have the following well-known relation:
$$J_L(t) \, = \, (-1)^{w(L)} t^{\frac{3w(L)-r(\wh G)+n(\wh G)}{4}}
    \,  \bigl(-t^{1/2}-t^{-1/2}\bigr)^{k(\wh G)-1}  \, T_{\Ga}(-t, -t^{-1}).
$$
\end{cor}

\bigskip
\section{Proof of the Main Theorem.}
The notation used in the definitions of the Kauffman bracket
(Definition~\ref{def:kb}) and the Bollob\'as-Riordan polynomial
(Definition~\ref{def:br}) hints on how to prove the
Main Theorem.
Since the crossings of the diagram $L$ correspond to the edges of $G=G_L$,
there is a
natural one-to-one correspondence $\vp: \cS(L) \to \cF(G)$
between the states $S \in \cS(L)$ and spanning
subgraphs $F\subseteq G$.
Namely,
let an $A$-splitting of a crossing in $S$ mean that we keep the
corresponding edge in the spanning subgraph $F= \vp(S)$.
Similarly, let a $B$-splitting in $S$ mean that we remove the
edge from the subgraph $F = \vp(S)$.

By definition, we have $\de(S) = \bc(F)$, for all $F = \vp(S)$.
Furthermore, we easily obtain the following
relation between the parameters:
$$e(F)=\a(S),\qquad e(G)-e(F)=\b(S)\ ,$$
for all $S \in \cS(L)$, and $F = \vp(S)$.
Now, for a spanning subgraph $F \in \cF(G)$, consider the term
$x^{r(G)-r(F)} y^{n(F)} z^{k(F)-\bc(F)+n(F)}$ of $R_G(x,y,z)$.
After a substitution
$$x\, = \, \frac{Bd}{A},\qquad y \, = \, \frac{Ad}{B},\qquad z\, =\,
\frac{1}{d}$$
and multiplication of this term by
$A^{r(G)} B^{n(G)} d^{k(G)-1}$ as in the Main Theorem,
we get
$$\begin{array}{ccl}
\multicolumn{3}{l}{ A^{r(G)} B^{n(G)} d^{k(G)-1}
                    (A^{-1}Bd)^{r(G)-r(F)} (AB^{-1}d)^{n(F)}
                    d^{-k(F)+\bc(F)-n(F)}  }     \vspace{12pt}\\
&=\!& A^{r(G)-r(G)+r(F)+n(F)} B^{n(G)+r(G)-r(F)-n(F)}
          d^{k(G)-1+r(G)-r(F)+n(F)-k(F)+\bc(F)-n(F)}  \vspace{12pt}\\
&=\!& A^{r(F)+n(F)} B^{n(G)+r(G)-r(F)-n(F)}
          d^{k(G)-1+r(G)-r(F)-k(F)+\bc(F)}\ .
\end{array}$$

It is easy to see that $r(F)+n(F)=e(F)$, and
$k(F)+r(F)=v(F)=v(G)$.  Therefore, $k(G)-k(F)+r(G)-r(F)=0$,
and we can rewrite our term as
$$A^{e(F)} B^{e(G)-e(F)} d^{\bc(F)-1}\ .$$
In terms of the state $S = \vp^{-1}(F) \in \cS(L)$
this term is equal to
$$A^{\a(S)} \, B^{\b(S)}\, d^{\de(S)-1}\,,$$
which is precisely the term of $\kb{L}$ corresponding to the
state $S \in \cS(L)$.  This completes the proof.
\hspace{\fill}$\square$

\bigskip

\section{Final remarks and open problems.}\label{rem_kb}
\mbox{\ }\vspace{-15pt}

{\bf 1.} \,
Trivalent ribbon graphs are the main objects in
the finite type invariant theory of knots, links and
3-manifolds, while general ribbon graphs appeared in the
literature under a variety of different names
(see e.g.~\cite{DKC,BR2,Ka1}).
Embeddings of ribbon graphs into 3-space are studied in~\cite{RT}.

\smallskip

{\bf 2.} \,
The Bollob\'as-Riordan polynomial can be defined by recurrent
contraction-deletion relations or by spanning tree expansion
similar to those of the Tutte polynomial, except that deletion of a loop
is not allowed.  We refer to~\cite{BR2,BR3} for the details.
We should note that~\cite{BR3} gives an extension to
unorientable surfaces as well. One can also find the contraction-deletion
relations and the spanning tree expansion for the signed
Bollob\'as-Riordan polynomial defined by (\ref{def_sbr}).
For a planar signed ribbon graph $\wh G$, the signed Bollob\'as-Riordan
polynomial $R_{\wh G}$ is related to
Kauffman's \emph{signed Tutte polynomial}~$Q[\wh G]$
from \cite{Ka2} by the formula
$$R_{\wh G}(x,y,z) \, = \, x^{\frac{v(\wh G)+1}{2}-k(\wh G)} \,
   y^{\frac{-v(\wh G)+1}{2}} \,
   Q[\wh G]\Bigl((y/x)^{1/2}, 1, (xy)^{1/2}\Bigr).
$$
So, our version of $R_{\wh G}$ may be considered as a generalization
of the polynomial~$Q[\wh G]$ to signed ribbon graphs.
If, besides planarity, all edges of $\wh G$ are positive, then
 $R_{\wh G}$ is related to the dichromatic polynomial $Z[\Ga](q,v)$
 (see \cite{Ka2}) of the underlying graph $\Ga$:
 $$R_{\wh G}(x,y,z) \, = \, x^{-k(\wh G)} \,  y^{-v(\wh G)} \,
 Z[\Ga](xy,y)\, .
 $$

\smallskip

{\bf 3.} \,
An interesting construction of a ribbon graph for a classical link
diagram is given in a recent paper~\cite{DFKLS}.  It allows to
compute the Jones polynomial of a general (not necessary alternating)
classical link using the ordinary (unsigned) version of the
Bollob\'as-Riordan polynomial.  
Most recently, Jeremy Voltz found a common generalization of~\cite{DFKLS}
and our Theorem~4.1 to arbitrary virtual links (in preparation).

\smallskip

{\bf 4.} \,
It would be interesting to generalize the Bollob\'as-Riordan polynomial
for colored ribbon graphs \cite{BR1,Tr} and prove the corresponding
relation with the Kauffman bracket.  Let us also mention that
in~\cite{Ja} (see also~\cite{Tr}), Jaeger found a different relation
between links and graphs and proved that the whole Tutte polynomial,
not just its specialization, can be obtained from the HOMFLY
polynomial of the appropriate link.  Recently Moffatt extended these results to
ribbon graphs \cite{Mof}.

Finally, recent results concerning combinatorial evaluations of
the Tutte and Bollob\'as-Riordan polynomials~\cite{KP} leave an
open problem of finding such evaluations for general values of the
polynomial~$R_G$.  It would be interesting to use the
Main Theorem to extend the results of~\cite{KP}.

\vskip.6cm

{\bf Acknowledgements}

\smallskip

\noindent
The first author thanks Alexander Shumakovich for pointing to
the reference \cite{Kam1}, and Oliver Dasbach, Xiao-Song Lin,
Jo Ellis-Monaghan, and  Neal Stoltzfus for useful conversations.
The second author is grateful to Jo Ellis-Monaghan,
Mike Korn and Vic Reiner for their insights into the
Tutte polynomial, and to B{\'e}la Bollob{\'a}s,
Vaughan Jones and Richard Stanley for the encouragement.
The second author was supported by the NSF.


\vskip1.cm

\parbox[t]{2.5in}{\it \textbf{Sergei~Chmutov}\\
Department of Mathematics\\
The Ohio State University, Mansfield\\
1680 University Drive\\
Mansfield, OH 44906\\
~\texttt{chmutov@math.ohio-state.edu}} \qquad\qquad
\parbox[t]{2.5in}{\it \textbf{Igor~Pak}\\
Department of Mathematics\\
Massachusetts Institute of Technology\\
77 Massachusetts Ave\\
Cambridge, MA 02139\\
~\texttt{pak@math.mit.edu}}

\end{document}